\newcommand{\F}{\mathcal{F}}
\newcommand{\B}{\mathcal{B}}
\newcommand{\Li}{\mathcal{L}}

\newcommand{\p}{\mathbb{P}}
\newcommand{\R}{\mathbb{R}}

\newcommand{\E}{\mathbb{E}}
\newcommand{\N}{\mathbb{N}}

\newcommand{\ds}{\displaystyle}

\newcommand{\li}{{\langle}}
\newcommand{\ri}{{\rangle}}
\def\1{\mbox{$1$\hskip-.23em\rule[1.33ex]{.17em}{.2ex}\hskip-.27em$1$}}

\documentclass[10pt,draft,reqno]{amsart}
     \makeatletter
     \def\section{\@startsection{section}{1}%
     \z@{.7\linespacing\@plus\linespacing}{.5\linespacing}%
     {\bfseries%\normalfont\scshape
     \centering
     }}
     \def\@secnumfont{\bfseries}
     \makeatother
\newtheorem{theorem}{Theorem}[section]
\newtheorem{lemma}[theorem]{Lemma}
\newtheorem{proposition}[theorem]{Proposition}

\theoremstyle{definition}
\newtheorem{definition}[theorem]{Definition}
\newtheorem{example}[theorem]{Example}
\theoremstyle{remark}
\newtheorem{remark}[theorem]{Remark}
\numberwithin{equation}{section}
\begin{document}

\title[Evolution systems of measures with L\'{e}vy noise]
{Evolution systems of measures for non-autonomous Ornstein-Uhlenbeck
processes with L\'{e}vy noise}

\author[Robert Wooster]{Robert Wooster*}
\thanks{* This research is partially supported by NSF Grant DMS-0706784 and SFB 701, Bielefeld University, Germany}
\address{Robert Wooster: Department of Mathematical Sciences, United States Military Academy, West Point, NY 10996, USA}
\email{robert.wooster@usma.edu}
\urladdr{http://www.dean.usma.edu/departments/math/people/wooster/}

%\author[Jane Doe]{Jane Doe*}
%\address{Jane Doe: Department, University, City, State with zip code, Country}
%\email{janedoe@dept.univ.edu}



\begin{abstract}
We begin by stating some of the important properties of L\'evy processes, and some well-known results about stochastic differential equations with L\'evy noise terms.  We then examine the question of existence and uniqueness of evolution
systems of measures for non-autonomous Ornstein-Uhlenbeck-type
processes with jumps. Finally we give some examples where we
explicitly compute the densities of such families of measures.
\end{abstract}

\maketitle

%\noindent $\clubsuit$ Note to author: Use 2000 Mathematics Subject Classification.

%\section{The First Section}

\section{Preliminaries}\label{s.1}

% Introduction

\subsection{Introduction}
In this paper we fix without further mention a filtered probability space
$(\Omega, \F, \F_{t}, \p)$ where the filtration $\F_{t}$ satisfies the usual hypothesis of completeness and right continuity. Consider the non-autonomous Ornstein-Uhlenbeck type
stochastic differential equation
% O-U SDE
\begin{align}\label{OU}
& dX(t) = (A(t) X(t-)+ f(t)) dt + B(t)dZ(t) \notag \\
& X(s) = x,
\end{align}
taking values in $\R^{d}$, and where $x \in \R^{d}$, $s \leq t \in \R$, $f: \R \to \R^{d}$,
$A: \R  \to\Li(\R^{d})$, $B: \R \to \Li(\R^{d})$ are bounded and continuous, and $Z: \Omega \times \R \to
\R^{d}$ is a L\'evy process. L\'evy processes are usually defined on $\R^{+}$.  We extend this to all of $\R$ by taking $(Z(t), t < 0)$ to be an independent copy of $-Z(-t-)$, see \cite{Applebaum09}, p. 239.  The process $Z$ will be referred to as the {\it noise} term in the equation.

A L\'evy process is a stochastic process with stationary and independent increments.  L\'evy processes are natural candidates to work with for two reasons.  First, examples are robust, including the stable (which include Gaussian) and Poisson processes.  The second reason is that the characteristic function has an elegant form, given by the celebrated L\'evy-Khintchine formula, Theorem \ref{LK}.  This allows us to study L\'evy processes, and the stochastic processes that solve (\ref{OU}), using Fourier analysis.

The goal of our paper is to give an existence and uniqueness result of an evolution system of
measures for the solution to (\ref{OU}).  The definition of
an evolution system of measures is given below and can be thought of
as a natural generalization of the notion of an invariant measure to
the non-autonomous case.  We pay particular attention to the case
where $Z$ is $\alpha$-stable with index of stability $0 < \alpha \leq 2$.  In the cases where $\alpha = 1$ and $\alpha=2$, we
explicitly compute the densities of evolution families.

Our paper is organized as follows.  In Section \ref{s.2} we give an explicit form for the solution of (\ref{OU}), which we denote $X_{s,x}(t)$.  The two-parameter transition evolution operator corresponding to $X_{s,x}(t)$ is defined as usual,
\begin{align*}
P_{s,t}f(x) := \E[f(X_{s,x}(t))] = \int_{\R^{d}} f(y) p_{s,t}(x,dy), \hspace{.1in} f \in B_{b}(\R^{d})
\end{align*}
where $p_{s,t}(x, A)$ is the transition probability of
$X_{s,x}(t)$, i.e. \[p_{s,t}(x,A) = \p(X(t) \in A | X(s) = x), A \in \B(\R^{d}).\]

The main result of this paper is in Section \ref{s.3}.  There we prove the existence of a unique family of probability measures, $\{ \nu_{t}\}_{t \in \R}$, which satisfy the equation
\[
\int_{\R^{d}}P_{s,t}f(x)\nu_{s}(dx) = \int_{\R^{d}}f(x)\nu_{t}(dx), \hspace{.1in} -\infty < s \leq t < \infty,
\]
for any $f \in B_{b}(\R^{d})$.  Such a family is called an {\it evolution family of measures}.
% F-R proof

Existence and uniqueness of invariant measures for processes that arise from solutions to autonomous versions of (\ref{OU}) are well-known.  For example the solution to the stochastic initial value problem,
\begin{align*}
dX(t) &= A X(t) + dZ(t) \\
X(0) &= x,
\end{align*}
is the Ornstein-Uhlenbeck process.

Here $A$ is the generator of a strongly continuous semigroup, $T_t$, of linear operators on a space $E$.  The transition semigroup of $X$ is given by {\it Mehler's formula}
\[
P_{t}f(x) = \int_{E}f(T_{t}x+y)\mu_t(dy), t \geq 0,
\]
where $\mu_t$ is a family of probability measures which satisfy
\[
\mu_{t+s} = (\mu_{t} \circ T_{s}^{-1}) \ast \mu_{s},\ \rm{for\ all}\ s,t \geq 0.
\]

Such processes were first studied in a Hilbert space setting by Chojnowska-Michalik in \cite{chojnowska87}.  Further work in this area was done by Applebaum in \cite{applebaum06} and \cite{applebaum07}, van Neerven \cite{neerven06}, and Lescot and Rockner \cite{lescot04}.  Existence of invariant measures for such processes are well-known.  In \cite{Rockner00}, Fuhrman and Rockner gave conditions under which a unique invariant measure exists for the semigroup $P_{t}$.

The techniques used in the proof of Theorem \ref{bigthm} are similar to the proof of Theorem 3.1 in \cite{Rockner00}.  In that paper Fuhrman and R\"ockner decomposed the law of $X$ into its deterministic, drift, and jump parts, and gave a proof based on weak compactness and weak convergence.  However because the result in this paper is for a non-autonomous setting, there are some significant differences.  One of them being a tightness condition on the evolution family of measures in Lemma \ref{lemma}.

Work done in the non-autonomous setting with the noise being a Gaussian process was carried out by DaPrato and Lunardi in \cite{DaPrato06}.  They proved that if $A(t)$ is $T$-periodic and $Z(t)$ is a $d$-dimensional Brownian
motion, then an evolution family of measures exists.  They also
showed that under the additional assumption that $\nu_{t}$ is
$T$-periodic, there exists a unique $T$-periodic evolution family of measures.
Because the noise term in \cite{DaPrato06} is a Brownian motion,
the evolution family of measures is Gaussian, and a formula for the
mean and variance is computed explicitly.  In this paper, besides considering a much more general class of driving noise, we do not make any
periodicity assumption, and our result agrees with DaPrato and Lunardi in this case.

In general, we cannot expect to be able to compute the laws of an
evolution family explicitly if we replace the noise term with a
L\'evy process other than a Brownian motion.  However, in Section
\ref{s.4} we give an example where we can explicitly compute
the laws of such an evolution system of measures if $Z$ is a symmetric $\alpha$-stable L\'evy
process, with $\alpha = 1$.  A value of $\alpha=2$ would yield a
result consistent with Da Prato and Lunardi, see \cite{DaPrato06}.  Even though we expect most of the results to be easily adjusted to the case when $\R^{d}$ is replaced by a separable Hilbert space, we restrict ourselves to the finite-dimensional case in part to give very explicit examples.

% notation

\subsection{Notation}
Throughout the paper we will use the following notation:

\vspace{.1in}

\begin{itemize}
\item $B_{r}(a) = \{ x \in R^{d}: |x - a| < r \} $ denotes the open ball of radius $r$ centered at $a \in \R^{d}$.  $\bar B_{R}$ is the closed ball.
\item $\mathcal{B}(\R^{d})$ denotes the Borel $\sigma$-field on $\R^{d}$.
\item $B_{b}(\R^{d})$ denotes the space of all bounded Borel functions on $\R^{d}$.
\item $C_{b}(\R^{d})$ denotes the space of all bounded continuous functions on $\R^{d}$.
\item $\Li(\R^{d})$ denotes the space of all $d \times d$ real-valued matrices.
\item $\mu \ast \nu$ denotes the convolution of two Borel probability measures on $\R^{d}$, \[(\mu \ast \nu)(A) = \int_{\R^d}\mu(A-x)\nu(dx),\]
for any $A \in \mathcal{B}(\R^{d}).$
\item $\hat \nu$ denotes the characteristic function of a probability measure $\nu$ on $\R^{d}$,
\[
\hat \nu(a) = \int_{\R^{d}}e^{i\langle a, x \rangle}\nu(dx).
\]
\end{itemize}

%\section{The Second Section}

\section{Ornstein-Uhlenbeck type stochastic differential equation}\label{s.2}

In this section we recall the notion of a strong solution to (\ref{OU}) and compute its characteristic function.

% Strong solution
\begin{definition}
A {\it strong solution} to (\ref{OU}) is a c\`adl\`ag process, adapted to the filtration generated by $Z(t)$, satisfying the integral equation,
\[
X(t) = x + \int_{s}^{t}(A(r)X(r-) + f(r))dr + \int_{s}^{t}B(r)dZ(r), \hspace{.2 in} s \leq t.
\]
We will sometimes denote this solution as $X_{s,x}(t)$.
\end{definition}

% Evolution operator associated to A(t)

Let $U(t,s)$ denote the {\it evolution operator} in $\R^{d}$ associated with $A(t)$. That is, $U(t,s)$ is the two parameter family of operators which solve
\begin{align*}
& \frac{\partial U(t,s)}{\partial t} = A(t) U(t,s),\ s \leq t, \notag \\
& U(s,s) = I,
\end{align*}
where $I$ is the identity operator.

It satisfies the properties
\begin{align*}
& U(t,s)U(s,r) = U(t,r), \hspace{.1in} r \leq s \leq t \in \R, \\
& U(s,r)^{T}U(t,s)^{T} = U(t,r)^{T}, \hspace{.1in} r \leq s \leq t \in \R, \\
& \frac{ \partial U(t,s)}{\partial s} = - U(t,s)A(s), \hspace{.1in} s \leq t \in \R,
\end{align*}

as shown in \cite{Pazy83}, p. 128-9.

We make the following stability assumption on evolution operator, $U(t,s)$.

% Assumption on U(t,s)

There exists $C, \epsilon > 0$ such that,
\begin{equation}\label{assumption1}
||U(t,s)|| \leq C e^{-\epsilon(t-s)},
\end{equation}
for all $-\infty < s \leq t < \infty$.

% Note on assumption 1

It is important to note that this cannot be replaced by an assumption on $A(t)$ itself.  For example, even if the eigenvalues of $A(t)$ are negative and bounded away from zero uniformly for all $t$, equation (\ref{assumption1}) need not hold, as can be found in \cite{Chi99} Example 3.5, p. 61.

% Solution to SDE theorem
\begin{theorem}
The stochastic differential equation \eqref{OU} has a unique strong solution for $t \geq s$, which we can write in terms of the evolution operator,
\begin{equation}\label{sol}
 X_{s,x}(t) =U(t,s)x + \int_{s}^{t}U(t,r)f(r)dr + \int_{s}^{t}U(t,r)B(r)dZ(r).
\end{equation}
\end{theorem}

\begin{proof}
Existence of a unique solution is given by Theorem 6.2.9, p. 374-5 of \cite{Applebaum09}.  Formula (\ref{sol}) is obtained by applying the variation of constants formula.
\end{proof}

%transition semigroup

The transition evolution operator is given by a generalized version of Mehler's formula
\[
P_{s,t}f(x) = \int_{\R^{d}}f(U(t,s)x+y)p_{s,t}(0,dy), \hspace{.1in} f \in B_{b}(\R^{d}).
\]

% L\'evy measure
\begin{definition}
A Borel measure, $M$, on $\R^{d}$ is called a {\it L\'evy measure} if
\[
\int_{\R^{d}}\frac{|y|^{2}}{1+|y|^{2}}M(dy) < \infty,
\]
and $M( \{0\} ) = 0$.  An equivalent definition sometimes used is
\[
\int_{\R^{d}}(1 \wedge |y|^2) M(dy) < \infty.
\]
\end{definition}

The next theorem gives the characteristic function of an infinitely divisible random variable in terms of three parameters.  In particular, for each fixed $t$ for a L\'evy process $Z$, $Z(t)$ has an infinitely divisible distribution.  The usefulness of the L\'evy-Khintchine formula is an important factor in choosing a L\'evy process for the noise term when working with stochastic differential equations with jumps.  A proof can be found, e.g. Theorem 1.2.14, p. 29 and Corollary 2.4.20, p. 127 in \cite{Applebaum09}.

% LK formula for L\'evy process
\begin{theorem}[L\'evy-Khintchine formula]\label{LK}
The characteristic function of the L\'evy process, $Z$, is of the form
\[
\phi_{Z(t)}(a) = \exp[-t\eta(a)],
\]
where
\[
\eta (a) = - i \li b,a \ri + \frac{1}{2} \li a, Ra \ri -
\int_{\R^{d}} \left [ e^{i\li a,y \ri}-1-\frac{i \li a,y
\ri}{1+|y|^{2}}  \right ]M(dy),
\]
$b \in \R^{d}$, R is a positive definite symmetric $d \times d$ matrix, and $M$ is a L\'evy measure on $\R^d - \{ 0 \}$.
The parameters $b,R,M$ are uniquely determined by the process and are called the {\it L\'evy triple} of $Z$.  Furthermore the function $\eta$ is continuous.

Conversely, any mapping of the form
\[
\phi(a) = \exp[-\eta(a)]
\]
is the characteristic function of an infinitely divisible random variable.
\end{theorem}

We now briefly pause to introduce the following notation.  We write $[b,R,M]$ to denote the probability law of an infinitely divisible random variable with triple $(b,R,M)$.  This is not standard but makes the notation in the proof of Theorem \ref{bigthm} easier.

% Characteristic function of solution

In the next two propositions we see that for each fixed $s,t \in \R$ and $x \in \R^{d}$, the solution to (\ref{OU}), $X_{s,x}(t)$, is an infinitely divisible random variable.  In Proposition \ref{prop1} we compute the characteristic function of this process.  The property that L\'evy processes have independent and stationary increments is important here.  In Proposition \ref{triple}, we utilize the L\'evy-Khintchine formula to compute the triple of $X_{s,x}(t)$.

\begin{proposition}\label{prop1}
The characteristic function of the process
\[
Y(t) = \int_{s}^{t}U(t,r)B(r)dZ(r)
\]
is of the form
\[ \phi_{Y(t)}(a) = \exp \left [-\int_{s}^{t}\eta(B(r)^{T}U(t,r)^{T}a)dr \right ]. \]
\end{proposition}

\begin{proof}
Fix $-\infty < s \leq t < \infty$.  Let $P_{n}=\{s = r_{0}^{(n)} \leq r_{1}^{(n)} \leq \cdot \cdot \cdot \leq r_{m(n)}^{(n)} = t\}$ be a sequence of partitions such that $||P_{n}|| \to 0$ as $n \to \infty$, where $||P_{n}||: = \max_{0\leq i \leq m(n)}(r^{(n)}_{i+1}-r^{(n)}_{i})$ is the mesh of the partition $P_{n}$.

By the construction of the It\^o stochastic integral,
\begin{align*}
\phi&_{Y(t)}(a)  = \E \exp \left ( i \langle a, Y(t) \rangle \right )\\
&= \E \exp \left ( i \langle a, \int_{s}^{t} U(t,r)B(r)dZ(r) \rangle \right ) \\
&= \E \exp \left ( i \left \langle a, \lim_{n \to \infty}
\sum_{j=1}^{m(n)} U(t,r_{j}^{(n)})B(r_{j}^{(n)})(Z(r_{j+1}^{(n)})-Z(r_{j}^{(n)}))
\right \rangle \right ).
\end{align*}
Next we take the limit out of the expectation using the Dominated Convergence theorem
\begin{align*}
\phi&_{Y(t)}(a) = \lim_{n \to \infty} \E \exp \left ( i \left \langle a, \sum_{j=1}^{m(n)} U(t,r_{j}^{(n)})B(r_{j}^{(n)})(Z(r_{j+1}^{(n)})-Z(r_{j}^{(n)})) \right \rangle \right )\\
& = \lim_{n \to \infty} \E  \prod_{j=1}^{m(n)} \exp  \left ( i \left
\langle a, U(t,r_{j}^{(n)})B(r_{j}^{(n)})(Z(r_{j+1}^{(n)})-Z(r_{j}^{(n)})) \right
\rangle \right ).
\end{align*}
In the next several steps we use the fact that $Z$ has independent
and stationary increments.
\begin{align*}
\phi_{Y(t)}&(a) = \lim_{n \to \infty} \prod_{j=1}^{m(n)} \E \exp \left ( i \left \langle a, U(t,r_{j}^{(n)})B(r_{j}^{(n)})(Z(r_{j+1}^{(n)})-Z(r_{j}^{(n)})) \right \rangle \right )\\
&= \lim_{n \to \infty} \prod_{j=1}^{m(n)} \E \exp \left ( i \left \langle B(r_{j}^{(n)})^{T}U(t,r_{j}^{(n)})^{T} a, (Z(r_{j+1}^{(n)})-Z(r_{j}^{(n)})) \right \rangle \right )\\
&= \lim_{n \to \infty} \prod_{j=1}^{m(n)} \E \exp \left ( i \left
\langle B(r_{j}^{(n)})^{T}U(t,r_{j}^{(n)})^{T} a, Z(r_{j+1}^{(n)}-r_{j}^{(n)}) \right
\rangle \right ).
\end{align*}
Finally we use Theorem \ref{LK} to finish the proof.
\begin{align*}
\phi_{Y(t)}(a) &=  \lim_{n \to \infty} \prod_{j=1}^{m(n)} \exp \left(- (r_{j+1}^{(n)}-r_{j}^{(n)}) \eta(B(r_{j}^{(n)})^{T}U(t,r_{j}^{(n)})^{T}a) \right )\\
& = \lim_{n \to \infty} \exp \left (- \sum_{j=1}^{m(n)} \eta(B(r_{j}^{(n)})^{T}U(t, r_{j}^{(n)})^{T}a)(r_{j+1}^{(n)}-r_{j}^{(n)}) \right ) \\
& =\exp \left (-  \lim_{n \to \infty} \sum_{j=1}^{m(n)} \eta(B(r_{j}^{(n)})^{T}U(t, r_{j}^{(n)})^{T}a)(r_{j+1}^{(n)}-r_{j}^{(n)}) \right ) \\
&= \exp \left (- \int_{s}^{t} \eta(B(r)^{T}U(t,r)^{T}a)dr \right ).
\end{align*}

\end{proof}

% Levy triple

\begin{proposition}\label{triple}
For each $ -\infty < s \leq t < \infty, x \in \R^{d}$, the random variable
$X_{s,x}(t)$ is infinitely divisible with the triple
\[ (U(t,s)x+b_{s,t}, R_{s,t}, M_{s,t}), \]
where
\begin{align*}
b_{s,t} &= \int_{s}^{t}U(t,r)f(r)dr + \int_{s}^{t}U(t,r)B(r)b\ dr \\
& + \int_{s}^{t} \int_{\R^{d}} U(t,r)B(r)y \left ( \frac{1}{1+|U(t,r)B(r)y|^{2}}-\frac{1}{1+|y|^{2}}  \right ) M(dy) dr,
\end{align*}
$
\ds R_{s,t} = \int_{s}^{t} U(t,r)B(r) R B(r)^{T}U(t,r)^{T}dr,\ {\rm and} \hfill
$
$
\ds M_{s,t}(A) = \int_{s}^{t}M(B(r)^{-1}U(t,r)^{-1}(A))dr.
$
\end{proposition}

\begin{proof}
Using Proposition \ref{prop1},
\begin{align*}
\phi&_{X_{s,x}(t)}(a) = \E \exp \left [ i \langle a,X_{s,x}(t) \rangle \right ]\\
=& \E \exp \left [ i \left \langle a,U(t,s)x + \int_{s}^{t} U(t,r)f(r)dr + \int_{s}^{t} U(t,r)B(r)dZ(r) \right \rangle \right ] \\
=& \exp \left [ i \left \langle a,U(t,s)x + \int_{s}^{t} U(t,r)f(r)dr \right \rangle -\int_{s}^{t}\eta(B(r)^{T}U(t,r)^{T}a)dr \right ].
\end{align*}
Now by Theorem \ref{LK},
\begin{align*}
\phi&_{X_{s,x}(t)}(a) = \exp \Biggl \{ i \left \langle a,U(t,s)x +  \int_{s}^{t} U(t,r)f(r)dr \right \rangle \\
& + \int_{s}^{t} \Biggl ( i \left \langle b,B(r)^{T}U(t,r)^{T}a \right \rangle \\
& - \frac{1}{2} \left \langle B(r)^{T}U(t,r)^{T}a, R B(r)^{T} U(t,r)^{T}a \right \rangle \\
& + \int_{\R^{d}} \left [ e^{i \left \langle B(r)^{T} U(t,r)^{T}a,y \right \rangle}-1-\frac{i \left \langle B(r)^{T} U(t,r)^{T}a,y \right \rangle}{1+|y|^{2}}  \right ]M(dy) \Biggl  )dr \Biggl \}.
\end{align*}
Next we rearrange some terms, and add and subtract $\frac{i \left \langle a,U(t,r)B(r) y \right \rangle}{1+|U(t,r)B(r) y|^{2}}$ to obtain
\begin{align*}
\phi&_{X_{s,x}(t)}(a) = \\
& \exp \Biggl \{ i \left \langle a,U(t,s)x +  \int_{s}^{t} U(t,r)f(r)dr + \int_{s}^{t}U(t,r) B(r) b \ dr \right \rangle \\
& - \frac{1}{2} \left \langle a, \int_{s}^{t} U(t,r)B(r) R B(r)^{T} U(t,r)^{T} dr \ a \right \rangle \\
& + \int_{s}^{t} \int_{\R^{d}} \biggl [ e^{i \left \langle a,U(t,r)B(r) y \right \rangle}-1-\frac{i \left \langle a,U(t,r)B(r) y \right \rangle}{1+|U(t,r)B(r) y|^{2}}\\
& + \frac{i \left \langle a,U(t,r)B(r) y \right \rangle}{1+|U(t,r)B(r) y|^{2}} - \frac{i \left \langle a,U(t,r)B(r) y \right \rangle}{1+|y|^{2}} \biggl ]M(dy) dr \Biggl \}.
\end{align*}
After more rearranging and a change of variables $\phi$ takes the desired form,
\begin{align*}
\phi&_{X_{s,x}(t)}(a) = \\
 & \exp \Biggl \{ i \left \langle a,U(t,s)x + \int_{s}^{t} U(t,r)f(r)dr + \int_{s}^{t}U(t,r) B(r) b \ dr \right \rangle \\
& + i \left \langle a , \int_{s}^{t} \int_{\R^{d}} U(t,r)B(r)y \left ( \frac{1}{1+|U(t,r)B(r) y|^{2}} - \frac{1}{1+|y|^{2}} \right ) M(dy) dr  \right \rangle \\
& - \frac{1}{2} \left \langle a, \int_{s}^{t} U(t,r)B(r) R B(r)^{T} U(t,r)^{T} dr \ a \right \rangle \\
& + \int_{s}^{t} \int_{\R^{d}} \biggl [ e^{i \left \langle a, z \right \rangle}-1-\frac{i \left \langle a, z \right \rangle}{1+|z|^{2}} M(B(r)^{-1}U(t,r)^{-1} dz)dr \Biggl \}. \\
\end{align*}
Now we show $R_{s,t}$ is non-negative definite, symmetric, and bounded, and that $M_{s,t}$ is a L\'evy measure.

Let $y \in \R^{d}$ and $s,t \in \R$.
\[
\left \langle y , R_{s,t} y \right \rangle = \int_{s}^{t} \left \langle B(r)^{T}U(t,r)^{T} y , R B(r)^{T}U(t,r)^{T} y \right \rangle dr \geq 0,
\]
since $R$ is non-negative definite.

Furthermore it easily follows that $R_{s,t}$ is symmetric since $R$ is symmetric.

%Furthermore since $R$ is symmetric,
%\begin{align*}
%\left \langle y , R_{s,t} y \right \rangle =& \int_{s}^{t} \left \langle y , U(t,r) B(r) R B(r)^{T}U(t,r)^{T} y \right \rangle dr \\
%=& \int_{s}^{t} \left \langle U(t,r) B(r) R B(r)^{T}U(t,r)^{T} y , y \right \rangle dr = \left \langle R_{s,t} y , y \right \rangle.
%\end{align*}
Let $|y| \leq 1$.  Since $B$ is bounded, let $C_{B}$ be such that $C_{B} \geq 1$ and $|| B(t) || \leq C_{B}$  for all $t \in \R$. By the stability assumption (\ref{assumption1}) we have,
\begin{align*}
| R_{s,t} y | &= \left | \int_{s}^{t} U(t,r)B(r)RB(r)^{T}U(t,r)^{T} dr \ y \right | \\
& \leq \int_{s}^{t} \left | U(t,r)B(r)RB(r)^{T}U(t,r)^{T} \ y \right | \ dr \\
& \leq C^{2} C_{B}^{2} ||R|| \int_{s}^{t} e^{-2 \epsilon(t-r)} dr = \frac{C^{2} C_{B}^{2} ||R||}{2 \epsilon} \left ( 1 - e^{-2 \epsilon(t-s)} \right ). \\
\end{align*}
Thus for all $-\infty < s \leq t < \infty$,
\[
|| R_{s,t} || = \sup_{|y| \leq 1} |R_{s,t}y| \leq \frac{C^{2} C_{B}^{2} ||R||}{2 \epsilon} < \infty.
\]
Since $M$ is L\'evy measure, set 
\[
K_{1} := \int_{ \{ |y| \leq 1 \}  } |y|^{2} M(dy),
\]
and
\[
K_{2} := M \left ( \left \{ |y| > \frac{1}{CC_{B}} \right \} \right ).
\]
Then for $ -\infty < s \leq t < \infty$, we have
\begin{align*}
\int_{\R^{d}} & \left ( 1 \wedge |y|^{2} \right ) M_{s,t}(dy) = \int_{s}^{t} \int_{\R^{d} } \left ( 1 \wedge |y|^{2} \right ) M(B(r)^{-1}U(t,r)^{-1} dy)dr \\
=& \int_{s}^{t} \int_{\R^{d} } \left ( 1 \wedge |U(t,r)B(r) z|^{2} \right ) M(dz)dr\\
\leq& \int_{s}^{t} \int_{\R^{d} } \left ( 1 \wedge C^{2}C_{B}^{2}|z|^{2} \right ) M(dz)dr \\
=& \int_{s}^{t} \biggl [ \int_{\left \{ |z| \leq \frac{1}{CC_{B}} \right \} } \left ( 1 \wedge C^{2}C_{B}^{2}|z|^{2} \right ) M(dz)\\
& + \int_{\left \{ |z| > \frac{1}{CC_{B}} \right \} } \left ( 1 \wedge C^{2}C_{B}^{2}|z|^{2} \right ) M(dz) \biggl ] dr \\
=&  (t-s) \left ( \int_{\left \{ |z| \leq \frac{1}{CC_{B}} \right \} } C^{2}C_{B}^{2}|y|^2 M(dy)dr + \int_{s}^{t} \int_{\left \{ |z| > \frac{1}{CC_{B}} \right \} }M(dy) \right ) \\
\leq & \left ( C^{2}C_{B}^{2} K_{1} + K_{2} \right )(t-s) < \infty,
\end{align*}
shows that $M_{s,t}$ are L\'evy measures.

This shows that the characteristic function of $X_{s,x}(t)$ is
\[
\phi_{X_{s,x}(t)}(a) = \exp(-\eta_{s,t,x}(a)),
\]
where,
\[
\eta_{s,t,x}(a)=\left ( -i \li b'_{s,t},a \ri + \frac{1}{2} \li a, R_{s,t} \ri - \int_{\R^{d}} \left [ e^{i\li a,y \ri}-1-\frac{i \li a,y \ri}{1+|y|^{2}}  \right ]M_{s,t}(dy) \right ),
\]
where $b'_{s,t}:=U(t,s)x + b_{s,t}$.  Thus $X_{s,x}(t)$ is infinitely divisible by the L\'evy-Khintchine formula.
\end{proof}

% triple as - \infty, t 

Note that as $s \to -\infty$, $M_{s,t}$ is an increasing family of measures.  Similarly $R_{s,t}$ is an increasing family of nonnegative symmetric matrices.  Because $B$ is bounded and by the stability assumption (\ref{assumption1}), $R_{s,t}$ converges as $s \to -\infty$.  We define
\[
R_{-\infty, t} := \int_{-\infty}^{t} U(t,r)B(r) R B(r)^{T} U(t,r)^{T} dr,
\]
and
\[
M_{-\infty,t}(A) := \sup_{s<t}M_{s,t}(A) = \int_{-\infty}^{t}M(B(r)^{-1}U(t,r)^{-1}(A))dr,
\]
$A \in \B(\R^{d})$.

% new section

\section{Invariant measures and evolution system of measures}\label{s.3}

Let us first consider the autonomous version of (\ref{OU})

\begin{align}\label{OUa}
& dY(t) = (A Y(t-)+z) dt + BdZ(t) \notag \\
& Y(s) = y,
\end{align}
 taking values in $\R^{d}$, and where $y,z \in \R^{d}$, $s,t \in \R$, $A,B \in \Li(\R^{d})$, and $Z: \Omega \times \R \to \R^{d}$ is a L\'evy process.

Let $P_{s,t}$ denote the transition evolution operator of a Markov process $X$
\[ P_{s,t}f(x) = \E [ f(X(t))|X(s) = x ] \] for $f \in B_{b}(\R^{d}).$

% time homogeneous

\begin{definition}
If $P_{s,t} = P_{0, t-s}$ for all $ 0 \leq s \leq t < \infty $, then $X$ is said to be {\it time-homogeneous}.  In this case we write $P_{0,t}$ as $P_{t}$.
\end{definition}
Since the coefficients $A, z,$ and $B$ in (\ref{OUa}) do not depend on $t$, the solution is a
time-homogeneous Markov process, as is shown in \cite{Applebaum09} Theorem 6.4.6, p. 388.

% Invariant measure
\begin{definition}
A Borel probability measure $\mu$ is an {\it invariant measure} for $X$ if
\[
\int_{\R^{d}}(P_{t}f)(x)\mu(dx) = \int_{\R^{d}}f(x)\mu(dx),
\]
for all $t \geq 0$, $f \in B_{b}(\R^{d})$.
\end{definition}

Unlike in the autonomous case, the solution to (\ref{OU}) is not time-homogeneous, so we cannot expect to find a single invariant measure.  Instead we look for a family of probability measures, ${\nu_{t}, t \in \R}$, called an {\it evolution family (or evolution system) of measures}.

% Evolution system of measures

\begin{definition}
Let $X_{s,x}(t)$ be as in (\ref{sol}).  A family of Borel probability measures, $\{ \nu_{t} \}_{t \in \R}$ is an {\it evolution family of measures} for the process, $X_{s,x}(t)$ if
\begin{equation}\label{evol}
\int_{\R^{d}}(P_{s,t}f)(x)\nu_{s}(dx) = \int_{\R^{d}}f(x)\nu_{t}(dx),
\end{equation}
for all $-\infty < s \leq t < \infty $, $f \in B_{b}(\R^{d})$.
\end{definition}

%remark

\begin{remark}\label{rem34}
Using a standard monotone class argument, (\ref{evol}) need only hold for indicator functions or for $f$ of the form $f(x) = \exp(i \li a,x \ri )$, $a \in \R^{d}$.
\end{remark}

% Lemma 2

\begin{lemma}\label{lemma2}
In the case of $f(x) = \exp(i \li a,x \ri )$, (\ref{evol}) becomes
\begin{align}\label{rem1}
\hat \nu_{s} \left (U(t,s)^{T} a \right )& \times \exp \left ( i \left \langle a , \int_{s}^{t} U(t,r)f(r)dr \right \rangle \right ) \times \notag \\
& \exp \left ( -\int_{s}^{t}\eta(B(r)^{T}U(t,r)^{T} a)dr \right ) = \hat \nu_{t}( a),
\end{align}
where $\hat \nu$ denotes the characteristic function of $\nu$.

In the case of $f = \1_{A}$, $A \in \mathcal{B}(\R^{d})$, (\ref{evol}) becomes
\begin{equation}\label{rem2}
\nu_{t}(A) = p_{s,t}(0,\cdot) \ast (\nu_{s} \circ U(t,s)^{-1} \cdot)(A).
\end{equation}
\end{lemma}

\begin{proof}
If $f(x) = \exp(i \li a,x \ri )$, then the left hand side of (\ref{evol}) becomes
\begin{align*}
\int_{\R^{d}} & (P_{s,t}f)(x)\nu_{s}(dx) =\int_{\R^{d}} \E \exp \left (i \li a,X_{s,x}(t) \ri \right )\nu_{s}(dx)\\
=& \int_{\R^{d}} \E \exp \Biggl (i \Bigl \langle a,U(t,s)x + \int_{s}^{t}U(t,r)f(r)dr \\
& + \int_{s}^{t}U(t,r)B(r)dZ(r) \Bigl \rangle \Biggl )\nu_{s}(dx).\\
\end{align*}
After writing the exponential of a sum as a product of exponentials, and pulling the deterministic factors outside the expectation, we obtain,
\begin{align*}
\int_{\R^{d}} (P_{s,t}f)&(x)\nu_{s}(dx) =\\
& \hat \nu_{s} \left (U(t,s)^{T} a \right ) \times  \exp \left ( i \left \langle a , \int_{s}^{t} U(t,r)f(r)dr \right \rangle \right ) \times \\
 & \E \exp \Biggl (i \Bigl \langle a, \int_{s}^{t}U(t,r)B(r)dZ(r) \Bigl \rangle \Biggl ).
\end{align*}
By Proposition \ref{prop1} the last term in the product is equal to
\[
 \exp \left (-\int_{s}^{t}\eta(B(r)^{T}U(t,r)^{T}a)dr \right ).
\]
This gives the left hand side of (\ref{rem1}).  The right hand side is obvious.

To prove (\ref{rem2}), let $f = \1_{A}$, for some $A \in \mathcal{B}(\R^{d})$.
The right hand side of (\ref{evol}) is then $\nu_{t}(A)$.  The left hand side of (\ref{evol}) works out to be
\begin{align*}
\int_{\R^{d}}(P_{s,t}f)(x)\nu_{s}(dx) &= \int_{\R^{d}} \int_{\R^{d}} \1_{A}(y)p_{s,t}(x,dy) \nu_{s}(dx)\\
&= \int_{\R^{d}} \int_{\R^{d}} \1_{A}(U(t,s)x + y)p_{s,t}(0,dy) \nu_{s}(dx)\\
&=  \int_{\R^{d}} \int_{\R^{d}} \1_{A}(x + y)p_{s,t}(0,dy) \nu_{s}(U(t,s)^{-1}dx)\\
&=  p_{s,t}(0,\cdot) \ast (\nu_{s} \circ U(t,s)^{-1} \cdot)(A).
\end{align*}
\end{proof}

% Lemma

\begin{lemma}\label{lemma}
Suppose the stability assumption (\ref{assumption1}) holds.  If $\{\nu_{s}\}_{s \in \R}$ is an evolution system
of measures for which there exists an integer $N_{0}$ such that the
sub-collection $\{\nu_{s}\}_{s<N_{0}}$ is uniformly tight, then
$\nu_{s} \circ U(t,s)^{-1} \to \delta_{0}$ weakly as $s \to
-\infty$ for each fixed $t$.
\end{lemma}

\begin{proof}
Fix $t$ and let $f \in C_{b}(\R^{d})$ and choose $M$ so that $||f|| \leq M$. Let
$\epsilon > 0$ be given.  Choose $\delta > 0$ so that if
$|x|<\delta$, then $|f(x)-f(0)|<\epsilon/2$.

Using the tightness assumption, choose $R$ so that
$$\ds \nu_{s}(\bar B_{R}(0)) > 1 - \frac{\epsilon}{4M}$$ for $s < N_{0}$.

Choose $N < N_{0} $, such that if $s < N$, then $\ds ||U(t,s)|| \leq
\frac{\delta}{R}$.  Then for any $x$ where $|x| \leq R$, we have that
$\ds |U(t,s)x| \leq ||U(t,s)|| \cdot |x| \leq \frac{\delta}{R}\ R = \delta.$

Then for $s < N$ we have
\begin{align}\label{lemeq}
&  \left | \int_{\R^{d}}  f(x) \nu_{s}(U(t,s)^{-1}dx)- \int_{\R^{d}}f(x) \delta_{0}(dx) \right | \\
&= \left | \int_{\R^{d}} f(U(t,s)x) \nu_{s}(dx)- \int_{\R^{d}}f(x) \delta_{0}(dx) \right | \notag \\
&= \left | \int_{\R^{d}} f(U(t,s)x) \nu_{s}(dx)- f(0) \right |. \notag
\end{align}
Since $\nu_{s}$ is a probability measure, (\ref{lemeq}) is equal to
\begin{align*}
& \left | \int_{\R^{d}} \left [ f(U(t,s)x) - f(0) \right ] \nu_{s}(dx) \right |  \\
& \leq \left | \int_{\bar B_{R}(0)}  \left [ f(U(t,s)x) - f(0) \right ] \nu_{s}(dx) \right | \\
& + \left | \int_{|x|>R}  \left [ f(U(t,s)x) - f(0) \right ] \nu_{s}(dx) \right | \\
& <  \epsilon/2 + 2M \cdot \epsilon/ 4M = \epsilon.
\end{align*}

\end{proof}

% Shift relative compactness and lemmas from Parthasarathy and Rockner

The next several results are taken from \cite{Pa67}, and are needed in the
proof of Theorem \ref{bigthm}.  We will also need the notion of
{\it shift relative compactness}.

% shift relative compactness

\begin{definition}
A set of Borel probability measures, $\mathcal H$ is said to be {\it shift relatively compact} if, for every sequence $\mu_{n} \in \mathcal H$, there is a sequence $\nu_{n}$ such that $\nu_{n}$ is a right (or left) translate of $\mu_{n}$, and $\nu_{n}$ has a convergent subsequence.
\end{definition}

% references from Parthasarathy

\begin{lemma}[Parthasarathy, Theorem III.2.2]\label{lem1}
Let $\{ \lambda_{n} \}, \{ \mu_{n} \}, \{ \nu_{n} \}$ be three sequences of measures on $\R^{d}$ such that $\lambda_{n} = \mu_{n} \ast \nu_{n}, n =1,2,... $If the sequence $\{ \lambda_{n} \}$ is relatively compact then the sequences $\{ \mu_{n} \}$ and $ \{ \nu_{n} \}$ are shift compact.
\end{lemma}

\begin{proof}
See p. 59, \cite{Pa67}.
\end{proof}

The following definition is the finite dimensional version of Definition VI.2.4, p. 155 in \cite{Pa67}.
\begin{definition}
A family $\{ S_{\alpha} \}$ of positive self-adjoint trace-class operators is said to be {\it compact} if $\sup_{\alpha} {\rm trace} (S_{\alpha})< \infty$.
\end{definition}

\begin{theorem}[Parthasarathy, Theorem VI.5.3]\label{thm1}
In order that a sequence $\mu_{n}$ of infinitely divisible distributions with representations $\mu_{n} = [ x_{n}, R_{n}, M_{n}]$ be relatively compact it is necessary and sufficient that the following hold:
\begin{itemize}
\item[(i)] $\{M_{n}\}$ restricted to to the complement of any neighborhood of the origin is weakly relatively compact.
\item[(ii)] $ \{S_{n} \}$ defined by
\begin{equation}\label{operator}
\li S_{n}y,y \ri = \li R_{-n,t}y, y \ri +\int_{|x| \leq 1}\li x,y \ri^{2}dM_{-n,t}(x)
\end{equation}
is compact.
\item[(iii)] $x_{n}$ is compact in $X$.
\end{itemize}
\end{theorem}

\begin{proof}
See p. 187, \cite{Pa67}.
\end{proof}

\begin{theorem}[Parthasarathy, Theorem III.2.1]\label{thm2}
Let $\{ \lambda_{n} \}, \{ \mu_{n} \}, \{ \nu_{n} \}$ be three sequences of measures on $\R^{d}$ such that $\lambda_{n} = \mu_{n} \ast \nu_{n}$ for each $n$.  If the sequences $\{\lambda_{n} \}$ and $\{ \mu_{n} \}$ are relatively compact then so is the sequence $ \{ \nu_{n} \}$.
\end{theorem}

\begin{proof}
See p. 58, \cite{Pa67}.
\end{proof}

We are now ready to prove the main result of this paper.  The result here generalizes Theorem 3.1 in \cite{Rockner00} to the time inhomogeneous case.

% Big theorem

\begin{theorem}\label{bigthm}
If there exists an evolution system of measures for $X_{s,x}(t)$ then the following conditions hold:
\begin{itemize}
\item[(i)] For any $\ds t \in \R, \sup_{s<t} {\rm tr} R_{s,t} < \infty,$
\item[(ii)] For any $\ds t \in  \R, \int_{-\infty}^{t}\int_{\R^{d}}(1 \wedge |U(t,r)B(r)y|^{2})M(dy)dr < \infty$.
\end{itemize}
If in addition,
\begin{itemize}
\item[(iii)] for any $t \in \R$, there exists an $N$ such that the collection $\{\nu_{t}\}_{t<N}$ is uniformly tight, then $\nu_t$ is unique and there exists
\[
b_{-\infty,t}:= \lim_{s \to -\infty} b_{s,t}.
\]
\end{itemize} 
Conversely if (i) and (ii) hold and $\lim_{s \to -\infty} b_{s,t}$ exists then for
each $t \in \R$, $M_{-\infty,t}$ is a L\'evy measure and there exists an
evolution system of measures, $\nu_{t}$, which is given by
\[
\nu_{t} = [b_{-\infty, t} , R_{-\infty, t}, M_{-\infty, t}].
\]
Recall from the paragraph following Theorem \ref{LK} that the bracket notation,
\[
[b_{-\infty, t} , R_{-\infty, t}, M_{-\infty, t}],
\]
denotes the law of infinitely divisible random variable with the L\'evy triple
\[
(b_{-\infty, t} , R_{-\infty, t}, M_{-\infty, t}).
\]
\end{theorem}

\begin{proof}

We prove the converse first.  Suppose (i), (ii) hold and the limit (iii) exists.
Fix $t \in \R$.  Using (ii),
\begin{align*}
\int_{\R^{d}}(1 \wedge |y|^{2})M_{-\infty,t}(dy) & = \int_{\R^{d}}(1 \wedge |y|^{2})\int_{-\infty}^{t}M(B(r)^{-1}U(t,r)^{-1}(dy))dr \\
& = \int_{-\infty}^{t} \int_{\R^{d}}(1 \wedge |y|^{2})M(B(r)^{-1}U(t,r)^{-1}(dy))dr \\
& = \int_{-\infty}^{t} \int_{\R^{d}}(1 \wedge |B(r)U(t,r)y|^{2})M(dy)dr < \infty
\end{align*}
shows that $M_{-\infty, t}$ is a L\'evy measure.

From the computation of the L\'evy triple of $X_{s,x}(t)$ in Proposition \ref{triple}, it follows that
$$ \hat \nu_{t}(a) = \exp \left ( i \left \langle a , \int_{-\infty}^{t} U(t,r)f(r)dr \right \rangle \right ) \exp \left \{-\int_{-\infty}^{t}\eta(B(r)^{T}U(t,r)^{T} a)dr \right \} .$$

Then using (\ref{rem1}) in Remark \ref{rem34},
\begin{align*}
\hat \nu_{s}  \biggl (U&(t,s)^{T} a \biggl ) \exp \left ( i \left \langle a , \int_{s}^{t} U(t,r)f(r)dr \right \rangle \right ) \times \\
& \exp \left \{-\int_{s}^{t}\eta(B(r)^{T}U(t,r)^{T} a)dr \right \} \\
=& \exp \left ( i \left \langle U(t,s)^{T} a , \int_{-\infty}^{s} U(s,r)f(r)dr \right \rangle \right ) \times \\
& \exp \left \{-\int_{-\infty}^{s}\eta(B(r)^{T}U(s,r)^{T} U(t,s)^{T} a)dr \right \} \times \\
&  \exp \left ( i \left \langle a , \int_{s}^{t} U(t,r)f(r)dr \right \rangle \right ) \exp \left \{-\int_{s}^{t}\eta(B(r)^{T}U(t,r)^{T} a)dr \right \} \\
=& \exp \left ( i \left \langle a , \int_{-\infty}^{t} U(t,r)f(r)dr \right \rangle \right ) \times \\
& \exp \left \{-\int_{-\infty}^{t}\eta(B(r)^{T}U(t,r)^{T} a)dr \right \}  = \hat \nu_{t}(a)
\end{align*}
shows that $\nu_{t}$ is an evolution system of measures.

Suppose now that an evolution system of measures, $\nu_{t}$, exists.
Fix $t$, then using  (\ref{rem2}) in Remark \ref{rem34}, for $s < t$,
\[
 \nu_{t} =
p_{s,t}(0, \cdot) \ast (\nu_{s} \circ U(t,s)^{-1}) = \delta_{b_{s,t}}
\ast [0,R_{s,t},0] \ast [0,0,M_{s,t}]  \ast (\nu_{s} \circ
U(t,s)^{-1}),
\]
where $\delta_{y}$ is the Dirac measure at $y$.

Set $s=-n$.  Then by Lemma \ref{lem1}, the sequence
$\delta_{b_{-n,t}} \ast [0,R_{-n,t},0] \ast [0,0,M_{-n,t}]$ is shift
relatively compact. This means that there is a sequence $y_{n} \in
\R^{d}$ (depending on $t$) such that
$$ \delta_{y_{n}} \ast \delta_{b_{-n,t}} \ast [0,R_{-n,t},0] \ast [0,0,M_{-n,t}] = [y_{n}+b_{-n,t},R_{-n,t},M_{-n,t}] $$
is weakly relatively compact.

Let $S_{n}: \R^{d} \to \R^{d}$ be a sequence of operators defined by (\ref{operator}).  By Theorem \ref{thm1}, the following hold:
\begin{itemize}
\item[(a)] $\{ M_{-n,t}\} $ restricted to the complement of any neighborhood of the origin is weakly relatively compact,
\item[(b)] $\sup_{n} {\rm tr}S_{n} <\infty$,
\item[(c)] $y_{n} + b_{-n,t}$ is relatively compact in $\R^{d}$.
\end{itemize}
%(see Definition VI.2.4 of \cite{Pa67}).
Part (a) implies
\[
M_{-\infty,t}(\{|x| \geq 1 \}) =
\sup_{n} M_{-n,t}(\{|x| \geq 1 \}) < \infty.
\]

By (b) we have
\begin{align*}
{\rm tr} R_{-\infty, t} + \int_{|x|\leq 1}|x|^{2}M_{-\infty,t}(dx) & = \sup_{n} \left ( {\rm tr}R_{-n,t} + \int _{|x|\leq 1}|x|^{2}M_{-n,t}(dx) \right ) \\
& = \sup_{n} {\rm tr} S_{n} < \infty.
\end{align*}
And so, by using Lemma 3.4 of \cite{Rockner00} for each fixed $t$, we have that $M_{-\infty,t}$ is a L\'evy measure and (i) and (ii) hold.

% converse

Now suppose also that there exists an $N$ such that the collection $\{\nu_{t}\}_{t<N}$ is uniformly tight.  Then by Lemma \ref{lemma}, $\nu_{s} \ \circ \ U(t,s)^{-1} \to \delta_{0}$ weakly as $s \to -\infty$.  By Lemma 3.4 \cite{Rockner00}, $[0,R_{s,t},0] \to [0,R_{-\infty,t},0]$ and $[0,0,M_{s,t}] \to [0,0,M_{-\infty,t}]$ weakly as $s \to -\infty$.  Thus by the weak continuity of convolution we conclude $$[0,R_{s,t},0] \ast [0,0,M_{s,t}] \ast (\nu_{s} \circ U(t,s)^{-1}) \to [0,R_{-\infty,t},0] \ast [0,0,M_{-\infty,t}].$$
Let $s_{n}$ be a sequence decreasing to $-\infty$.  Then
$$\nu_{t} = \delta_{b_{s,t}} \ast [0,R_{s,t},0] \ast [0,0,M_{s,t}]  \ast (\nu_{s} \circ U(t,s)^{-1}), $$
and by Theorem \ref{thm2}, the collection $\{ \delta_{b_{s_{n},t}} \}_{n \in \N}$ is weakly relatively compact.
Thus there is a probability measure $\sigma_{t}$ and a subsequence $n_{k}$ such that $\delta_{b_{s_{n_{k}},t}} \to \sigma_t $ weakly.  Letting $k \to \infty$,
$$ \nu_{t} = \sigma_{t} \ast [0,R_{-\infty,t},0] \ast [0,0,M_{-\infty,t}]. $$
Taking Fourier transforms of both sides we have
$$\hat \sigma_{t} = \hat \nu_{t} ( \widehat{[0, R_{-\infty,t},0 ]}  \cdot \widehat{ [0,0,M_{-\infty,t}]})^{-1},$$
We see that $\sigma_{t}$ does not depend on the subsequence, and
so $\delta_{b_{s_{n},t}}$ converges weakly.  This implies that
$b_{-\infty,t} := \lim_{n \to \infty} b_{s_{n},t}$ exists.  Since
$s_{n}$ is arbitrary, we have $b_{-\infty,t} = \lim_{s \to
-\infty}b_{s,t}$,

Thus  we have shown that \[ \nu_{t} = \delta_{b_{-\infty, t}} \ast
[0,R_{-\infty,t},0] \ast [0,0,M_{-\infty,t}] \] is uniquely
determined.

\end{proof}

% new section

\section{Examples}\label{s.4}

%ref examples

In this final section we give some examples where we explicitly compute the characteristic functions and densities of evolution systems of measures to which Theorem \ref{bigthm} applies.   The first example is where the noise term is a $d$-dimensional Gaussian process.  If we further require the coefficients to be $T$-periodic and take the Gaussian process to a be Brownian motion, the result agrees with  DaPrato and Lunardi in \cite{DaPrato06}.  In the second example we consider the case where $Z(t)$ is a one dimensional symmetric $\alpha$-stable process.  Where $Z(t)$ has a Cauchy distribution, i.e. $\alpha =1$, we explicitly compute the densities of the evolution system of measures.

% Gaussian example

\begin{example} \label{ex1}
Let $Z(t)$ have the L\'evy triple $(b,R,0)$, i.e. $Z(t)$ is a $d$-dimensional Gaussian process with mean vector $b$ and covariance matrix $R$.

The computation in Proposition \ref{triple} shows that the collection $\{ R_{s,t}, -\infty < s \leq t < \infty \}$ are non-negative definite, symmetric and uniformly bounded.  Thus $\sup_{s<t}{\rm tr} R_{s,t} < \infty$ and (\ref{assumption1}) implies that $\lim_{s \to -\infty} b_{s,t}$ exists.  Thus by Theorem \ref{bigthm}, the family of Gaussian measures with triple $\nu_{t} \sim [b_{-\infty,t}, R_{-\infty,t}, 0]$ is an evolution family of measures for the process $X_{s,x}(t)$, where
\begin{align*}
b_{-\infty,t} &= \int_{-\infty}^{t}U(t,r)f(r)dr + \int_{-\infty}^{t}U(t,r)B(r)b\ dr,\\
R_{-\infty,t} &= \int_{-\infty}^{t} U(t,r)B(r) R B(r)^{T}U(t,r)^{T}dr.\\
\end{align*}
\end{example}

% prop

The fact that stochastic integrals with respect to symmetric $\alpha$-stable processes are $\alpha$-stable make them a very useful subclass of L\'evy processes.  The next proposition summarizes this result.  A proof can be found in Samorodnitsky and Taqqu's book on stable processes, \cite{Sam94}.

% char function of solution (alpha stable)

\begin{proposition}\label{sam}
Suppose $Z(t)$ is a symmetric, $\alpha$-stable L\'evy process with characteristic function
\[
\phi_{Z(t)}(a) = e^{-\sigma^{\alpha}|a|^{\alpha}}.
\]
Then for each fixed $s \leq t$, the random variable
\[
Y_{s}(t) = \int_{s}^{t}e^{-\int_{u}^{t}\lambda(r)dr}dZ(u)
\]
has an $\alpha$-stable distribution with characteristic function
\[
\phi_{Y_{s}(t)}(a) = \exp \left \{ \sigma^{\alpha} \left ( \int_{s}^{t}e^{-\alpha \int_{u}^{t}\lambda(r)dr}du \right ) |a|^{\alpha} \right \}.
\]

\end{proposition}

%Stochastic differential equation equations with a symmetric $\alpha$-stable noise term were studied by Zanzotto \cite{zanzotto02}.

% alpha-stable

\begin{example} \label{ex2}
Let $Z(t)$ be a 1-dimensional symmetric $\alpha$-stable process with index of stability $0 < \alpha < 2$.

A one dimensional version of (\ref{OU}) is
\begin{align}\label{touex}
dX(t) &= \lambda (t) \left [ \mu(t) - X(t-) \right ] dt + \sigma (t) dZ(t),\notag \\
X(s) &= x,
\end{align}
where $\lambda, \mu, \sigma$ are bounded and continuous on $\R$, and $x \in \R$.  In addition we require $\lambda(t) \geq \epsilon > 0$ for all $t \in \R$.  Here the evolution operator has the form $U(t,s) = e^{-\int_{s}^{t}\lambda (r)dr}$.  The positivity condition on $\lambda$ implies that the stability assumption (\ref{assumption1}) is satisfied.

We write the solution to (\ref{touex}),
\begin{align*}
X(t) = X_{s,x}(t) =& e^{ -{\int_{s}^{t}\lambda (u)du}}x + \int_{s}^{t} e^{ -{\int_{r}^{t}\lambda (u)du}} \lambda(r)\mu(r)dr \notag \\
&+ \int_{s}^{t} e^{ -{\int_{r}^{t}\lambda (u)du}}\sigma(r)dZ(r).
\end{align*}
The transition evolution operator associated with $X$ takes the form,
\[
P_{s,t}f(x) = \E f(X_{s,x}(t)) = \int_{-\infty}^{\infty} f \left ( e^{-\int_{s}^{t} \lambda(r) dr} x + y \right ) p_{s,t}(0,dy).
\]
For each $t > s$, the law of $X_{s,x}(t)$ is $\alpha$-stable and has characteristic function
\begin{align*}
\phi_{X_{s,x}(t)} (a) &=  \exp \biggl [ i \left (e^{ -{\int_{s}^{t}\lambda (u)du}}x + \int_{s}^{t} e^{- \int_{r}^{t}\lambda (u)du} \lambda(r)\mu(r)dr \right ) a  \\
 &- \int_{s}^{t} e^{ - \alpha \int_{r}^{t}\lambda (u)du}[\sigma(r)]^{\alpha} dr \ |a|^{\alpha} \biggl ].
 \end{align*}
 The collection of measures, $\{ \nu_{t}, t \in \R \}$, with characteristic functions
 \[
 \hat \nu_{t} =  \exp \biggl [ i \int_{-\infty}^{t} e^{- \int_{r}^{t}\lambda (u)du} \lambda(r)\mu(r)dr \ a - \int_{-\infty}^{t} e^{ - \alpha \int_{r}^{t}\lambda (u)du}[\sigma(r)]^{\alpha} dr \ |a|^{\alpha} \biggl ],
 \]
 is the unique evolution system of measures for $X_{s,x}(t)$.
 In the case where $\alpha = 1$, we can explicitly write the densities of the $\nu_{t}$,
\[
	 f_{\nu_{t}}(y) =  \frac{a(t)}{ \pi \left [ \left (y- b(t) \right )^2+ \left (a(t) \right )^2 \right ],  }
	\]
	where 
	\[
	a(t) = \int_{-\infty}^{t}e^{ -\int_{r}^{t}\lambda (u)du} \sigma(r) dr,
	\]
	and
	\[
	b(t) = \int_{-\infty}^{t} e^{ -\int_{r}^{t}\lambda (u)du} \lambda(r) \mu(r) dr.
	\]
\end{example}

%%%%%%%%%%%

\par\bigskip\noindent
{\bf Acknowledgment.} I would like to thank Professor M. R\"ockner for many
valuable discussions during my stay at Bielefeld University in
summer of 2008 and Professor M. Gordina for introducing me to this exciting field and providing valuable insight and guidance.  Finally a special thanks also goes to the reviewer, whose comments greatly improved the exposition of this paper.

\bibliographystyle{amsplain}

\begin{thebibliography}{99}

\bibitem{applebaum06} Applebaum, D.:
{\em Martingale-valued Measures, Ornstein-Uhlenbeck Processes with Jumps and Operator Self-decomposability in Hilbert Space, In memoriam Paul-Andr\'e Meyer: S\'eminaire de Probabilit\'es XXXIX, Lecture Notes in Math.}, {\bf 1874}, Springer, Berlin, (2006) 171--196.

\bibitem{applebaum07} Applebaum, D.:
On the Infinitesimal Generators of Ornstein-Uhlenbeck Processes with Jumps in Hilbert Space, {\em Potential Anal.} {\bf 26} (2007) 79--100.

\bibitem{Applebaum09} Applebaum, D.:
{\em L\'evy Processes and Stochastic Calculus}, Cambridge University Press, Cambridge, 2009.

\bibitem{Chi99} Chicone, C., Latushkin, Y.:
{\em Evolution Semigroups in Dynamical System and Differential Equations}, American Mathematical Society, Providence, RI, 1999.

\bibitem{chojnowska87} Chojnowska-Michalik, A.:
On Processes of Ornstein-Uhlenbeck type in Hilbert Space, {\em Stochastics} {\bf 21} (1987) 251--286.

\bibitem{DaPrato06} Da Prato, G., Lunardi, A.:
Ornstein-Uhlenbeck Operators with Time Periodic Coefficients, {J. Evol. Equ.} {\bf 7} (2007) 587--614.

\bibitem{Rockner00} Fuhrman, M., R\"ockner, M.:
Generalized Mehler Semigroups: the Non-Gaussian Case, {\em Potential Anal.} {\bf 12} (2000) 1--47.

\bibitem{lescot04} Lescot, P., R\"ockner, M.:
Perturbations of Generalized Mehler Semigroups and Applications to Stochastic Heat Equations with L\'evy Noise and Singular Drift, {\em Potential Anal.} {\bf 20} (2004) 317--344.

\bibitem{Pa67} Parthasarathy, K.R.:
{\em Probability Measures in Metric Spaces}, Academic Press, New York, 1967.

\bibitem{Pazy83} Pazy, A,:
{\em Semigroups of Linear Operators and Applications to Partial Differential Equations}, Springer-Verlag, New York, 1983.

\bibitem{Sam94} Samorodnitsky, G., Taqqu, M.:
{\em Stable Non-Gaussian Random Processes}, Chapman \& Hall, New York, 1994.

\bibitem{neerven06} van Neerven, J. M. A. M., Weis, L.:
Invariant Measures for the Linear Stochastic Cauchy Problem and $\R$-boundedness of the Resolvent, {\em J. Evol. Equ.} {\bf 6} (2006) 205--228.

\end{thebibliography}

\end{document}